\documentclass[11pt,twoside]{article}
\usepackage[left=3cm,right=3cm,top=3cm,bottom=3cm]{geometry}
\usepackage{amssymb}
\usepackage{graphicx}

\usepackage{graphicx,type1cm,eso-pic,color}
\usepackage{amssymb}
\usepackage{amssymb,amsmath}
\usepackage{graphicx,epsfig}
\usepackage{enumerate}
\usepackage{color}
\usepackage {multicol}

\numberwithin{equation}{section}
\newtheorem{theorem}{Theorem}

\newtheorem{proposition}{Proposition}
\newtheorem{remark}{Remark}

\newtheorem{example}{Example}

\makeatletter
\AddToShipoutPicture{
	\setlength{\@tempdimb}{.5\paperwidth}
	\setlength{\@tempdimc}{.5\paperheight}
	\setlength{\unitlength}{1pt}
	\put(\strip@pt\@tempdimb,\strip@pt\@tempdimc)
}
\makeatother

\usepackage{scalerel,stackengine}
\stackMath
\newcommand\reallywidehat[1]{%
	\savestack{\tmpbox}{\stretchto{%
			\scaleto{%
				\scalerel*[\widthof{\ensuremath{#1}}]{\kern-.6pt\bigwedge\kern-.6pt}%
				{\rule[-\textheight/2]{1ex}{\textheight}}
			}{\textheight}%
		}{0.5ex}}%
	\stackon[1pt]{#1}{\tmpbox}%
}
\usepackage{mathtools,leftindex,tensor,mhchem}

\begin{document}
	\setcounter{page}{1}

	\thispagestyle{empty}
	\markboth{}{}

	\pagestyle{myheadings}
	\markboth{ S.K.Chaudhary }{ S.K.Chaudhary}
	
	\date{}
	
	
	\noindent  
	
	\vspace{.1in}

	{\baselineskip 20truept
		
		\begin{center}
			{\Large {\bf Study of inaccuracy measures of record values}} \footnote{\noindent	{\bf 1.}  Corresponding author E-mail: skchaudhary1994@gmail.com, skchaudhary1994@kgpian.iitkgp.ac.in, skc@bitmesra.ac.in\\
				{\bf 2. } E-mail: nitin.gupta@maths.iitkgp.ac.in}
			
	\end{center}}

	\vspace{.1in}
	
	\begin{center}
		{\large {\bf Santosh Kumar Chaudhary$^1$, Nitin Gupta$^2$}}\\
		
		{\large {\it ${}^{1}$	Centre for Quantitative Economics and Data Science, Birla Institute of Technology Mesra, Ranchi, Jharkhand 835215, India. }}\\	
		{\large {\it ${}^{2}$ Department of Mathematics, Indian Institute of Technology Kharagpur, West Bengal 721302, India. }} \\
	\end{center}

	\vspace{.1in}
	\baselineskip 12truept

	\begin{abstract}
		In this paper, we investigate inaccuracy measures based on record values, focusing on the relationship between the distribution of the $n$-th upper and lower $k$-record values and the parent distribution. We extend the classical Kerridge inaccuracy measure, originally developed for comparing two distributions, to record values and derive expressions for both upper and lower record values. In addition, we explore various other inaccuracy-based measures, such as cumulative residual inaccuracy, cumulative past inaccuracy, and extropy inaccuracy measures, and their applications in characterizing symmetric distributions. We compute these measures through illustrative examples for several well-known lifetime distributions, including the exponential, Pareto, and Weibull distributions. Our findings provide insights into how inaccuracy varies with record order and distribution parameters, contributing to a deeper understanding of information-theoretic measures applied to records.
		
		\vspace{.1in}
		
		
	\end{abstract}
	
	\section{Introduction}
	Inaccuracy measures have been instrumental in information theory and statistics for quantifying the deviation or discrepancy between different probability distributions. These measures help understand how much information is lost or gained when moving from one distribution to another. One of the foundational inaccuracy measures is Shannon entropy, which quantifies uncertainty and has paved the way for other important measures like Kerridge inaccuracy. These inaccuracy measures have applications in areas such as coding theory, statistical inference, and reliability analysis. This review explores the specific inaccuracy measures related to record values, including Kerridge inaccuracy, extropy, and cumulative variants.

	Inaccuracy measures, including Kerridge’s measure and Shannon’s entropy, evaluate the accuracy of statistical models. These measures have been comprehensively examined in coding theory by Nath (1968) and Bhatia (1995). Recent research on inaccuracy measures is reflected in the works of Kundu and Nanda (2015), Kundu et al. (2016), Kayal et al. (2017), and Psarrakos and Di Crescenzo (2018). Furthermore, Kumar and Taneja (2015) proposed a cumulative residual inaccuracy measure.    Tahmasebi et al. (2018) proposed a cumulative measure of inaccuracy in lower record values and study characterization results in the case of dynamic cumulative inaccuracy. Thapliyal and Taneja (2013) considered a measure of inaccuracy between distributions of the $i$-th order statistics and parent random variable. Mohammadi and Hashempour (2024) introduced the concept of dynamic cumulative residual extropy inaccuracy (DCREI) by expanding on the existing dynamic cumulative residual extropy (DCRE) measure and proposing a weighted version of it.
	
	Let $X$ and $Y$ be two non-negative random variables with cumulative distribution functions (cdf) $F_X(x),$ $F_Y(x),$ the probability density function(pdf) $f_X(x),$ $f_Y(x),$ and
	reliability functions $\bar{F}_X(x),$ $\bar{F}_Y(x),$ respectively.
	Kerridge (1961) introduced a measure of inaccuracy called Kerridge inaccuracy, which is a generalization of the Shannon entropy measure. If $F_X(x)$ is the actual distribution function corresponding to the observations and $F_Y(x)$ is the distribution assigned by the experimenter, the Kerridge inaccuracy measure of $X$ and $Y$ is defined as 
	\begin{align}\label{1}
		H(f_X,f_Y)=- \int f_X(x)\log f_Y(x)dx.
	\end{align}
	We point out that $H(f_X,f_X)=H(X)=- \int f_X(x)\log f_X(x)dx$ is entropy defined by Shannon (1948). Ahmadi (2021) characterized symmetric distribution using Kerridge inaccuthe racy measure of record values.

	The inaccuracy extropy measure (KIJ) is defined as
	\begin{align}\label{2}
		J(X,Y)=-\frac{1}{2} \int f_X(x)f_Y(x)dx.
	\end{align}
	We point out that $J(X,X)=J(X)=-\frac{1}{2} \int f^2_X(x)dx$  is extropy defined by Lad et al. (2015). Gupta and Chaudhary (2023) provided a characterization of symmetric distribution using inaccuracy extropy measure of record values.
	
	Cumulative residual extropy inaccuracy (CRIJ) and cumulative past extropy inaccuracy (CPIJ) measures, which are extensions of the corresponding cumulative residual extropy and Cumulative past extropy, respectively, are defined as :
	\begin{align}
		CRIJ(X,Y)&=-\frac{1}{2} \int \bar{F}_{X} (x)\bar{F}_Y(x)dx, \\
		\text{and} \ \	
		CPIJ(X,Y)&=-\frac{1}{2} \int F_{X} (x)F_Y(x)dx.
	\end{align}
	We point out that $CRIJ(X,X)$ and $CPIJ(X,X)$ becomes cumulative past extropy(CPJ) and  cumulative residual extropy(CRJ) respectively, given as
	\begin{align}
		CRIJ(X,X)=CPJ(X)&=-\frac{1}{2} \int \bar{F}^2_{X} (x)dx, \\
		\text{and} \ \	
		CPIJ(X,Y)=CRJ(X)&=-\frac{1}{2} \int F_{X}^2 (x)dx.
	\end{align}

	Gupta and Chaudhary (2023) provided a characterization of symmetric distribution using CPIJ and CRIJ of record values.
	
	In information theory, the Kullback-Leibler discrimination information measure proposed by Kullback and Leibler (1951) is a significant measure that has been used in reliability analysis and other related fields. It measures how closely (similarly) two statistical distributions resemble one another. The Kullback-Leibler divergence between  random variable $X$  and random variable $Y$ as
	\begin{align}\label{defkullback}
		K(X,Y)=\int f_X(x)\ln\left(\frac{f_X(x)}{f_Y(x)}\right)dx. 
	\end{align}

	The amount of relative information between  random variable $X$  and random variable $Y$ as is determined by

	\begin{align}\label{3}
		I(X,Y)=\frac{1}{2} \int f_X(x) \left(f_X(x)-f_Y(x)\right)dx.
	\end{align}
	
	Note that $K(X,X)=0$ and $I(X,X)=0.$

	Kumar and Taneja (2012) proposed a cumulative residual inaccuracy measure as 
	\begin{align}\label{4}
		H(\bar{F}_X,\bar{F}_Y)=- \int \bar{F}_X(x)\log \bar{F}_Y(x)dx.
	\end{align}
	
	Thapliyal and Taneja (2015) proposed a cumulative past inaccuracy measure as 
	\begin{align}\label{3}
		H({F}_X,{F}_Y)=- \int F_X(x)\log F_Y(x)dx.
	\end{align}

	We point out that $\xi(X,X)=\xi(X)=- \int \bar{F}_X(x)\log \bar{F}_X(x)dx$ and  $\bar \xi(X,X)=\bar \xi(X)=- \int F_X(x)\log F_X(x)dx$ are cumulative residual entropy defined by  Rao et al. (2004) and cumulative past entropy defined by Di Crescenzoa and Longobardi (2009), respectively.

	\indent The concept of $k$-records was introduced by Dziubdziela and Kopociński (1976); for more details, also see  Ahsanullah (1995) and Arnold et al. (1998). The pdf of the $n$th upper record value $U_{n}$ and the $n$th lower record value $L_{n}$ respectively are given by (see Arnold et al. (2008) and  Ahsanullah (2004)) 
	\begin{align}
		&f_{U_{n}}(x)=\dfrac{1}{(n-1)!}[-\log \bar{F}_X(x)]^{n-1} f_X(x),\  -\infty<x<\infty, \label{funk}\\
		\text{and} \ \ \ &f_{L_{{n}}}(x)=\dfrac{1}{(n-1)!}[-\log F_X(x)]^{n-1} f_X(x),\  -\infty<x<\infty  \label{flnk}.
	\end{align}
	\noindent The cdf of $U_{n}$ and  $L_{n}$, respectively, are  
	\begin{align}
		F_{U_{n}}(x)&=  1-\bar{F}_X (x)   \sum_{i=0}^{n-1} \frac{(- \log\bar{F}_X (x) )^i}{i!}, \label{FUnku}\\
		\text{and} \ \ \ 	F_{L_{n}}(x)&=  {F}_X (x)   \sum_{i=0}^{n-1} \frac{(- \log F_X (x) )^i}{i!}.\label{FLnku}
	\end{align}
	The pdf of the $n$th upper $k$-record value $U_{n,k}$ and the $n$th lower $k$-record value $L_{n,k}$ respectively are given by (see Arnold et al. (2008) and  Ahsanullah (2004)) 
	\begin{align}
		&f_{U_{n,k}}(x)=\dfrac{k^n}{(n-1)!}[-\log \bar{F}_X(x)]^{n-1} [\bar{F}_X(x)]^{k-1}f_X(x),\  -\infty<x<\infty, \label{funk}\\
		\text{and} \ \ \ &f_{L_{{n,k}}}(x)=\dfrac{k^n}{(n-1)!}[-\log F_X(x)]^{n-1} [F_X(x)]^{k-1}f_X(x),\  -\infty<x<\infty  \label{flnk}.
	\end{align}
	\noindent The cdf of $U_{n,k}$ and  $L_{n,k}$, respectively, are  
	\begin{align}
		F_{U_{n,k}}(x)&=  1-\bar{F}^k_X (x)   \sum_{i=0}^{n-1} \frac{(-k \log\bar{F}_X (x) )^i}{i!}, \label{FUnku}\\
		\text{and} \ \ \ 	F_{L_{n,k}}(x)&=  {F}^k_X (x)   \sum_{i=0}^{n-1} \frac{(-k \log F_X (x) )^i}{i!}.\label{FLnku}
	\end{align}

	\section{Kerridge inaccuracy measure of record values}
	Many authors studied the measure of inaccuracy between distributions of the records as well as order statistics and parent random variable (See, Thapliyal and Taneja (2015) and Goel et al. (2018a, (2018b)). The Kerridge inaccuracy measure associated with $n$th upper $k$-record value $U_{n,k}$  and parent density function $f(x)$ is given as
	\begin{align}\label{1}
		H(f_{U_{n,k}},f_X)&=- \int_{-\infty}^{\infty} f_{U_{n,k}}(x)\log f_X(x)dx \nonumber  \\
		&=-\int_{-\infty}^{\infty} \dfrac{k^n}{(n-1)!}[-\log \bar{F}_X(x)]^{n-1} [\bar{F}_X(x)]^{k-1}f_X(x) \log f_X(x)dx	\nonumber\\
		&= -\int_{0}^{\infty} \dfrac{k^n}{(n-1)!}t^{n-1} e^{-kt} \log f_X(F^{-1}(1-e^{-t}))dt	      \nonumber  \\
		&=E(-\log f_X(F_X^{-1}(1-e^{-T_{n,k}})))		         
	\end{align}
	where $T_{n,k}$ has Gamma distribution with parameters $n$ and $k$ having pdf 
	\begin{equation*}
		f_{T_{n,k}}(t)=\dfrac{k^n}{(n-1)!}t^{n-1} e^{-kt}, \ t>0, \ n>0, \ k>0.
	\end{equation*}
	Here, $n$ and $k$ have positive integer values.
	The Kerridge inaccuracy measure associated with $n$th lower $k$-record value $L_{n,k}$  and parent density function f(x) is given as
	\begin{align}\label{2}
		H(f_{L_{n,k}},f_X)&=- \int_{-\infty}^{\infty} f_{L_{n,k}}(x)\log f_X(x)dx \nonumber  \\
		&=- \int_{-\infty}^{\infty} \dfrac{k^n}{(n-1)!}[-\log F_X(x)]^{n-1} F_X(x)^{k-1}f_X(x) \log f_X(x)dx	\nonumber\\
		&=-	\int_{0}^{\infty} \dfrac{k^n}{(n-1)!}t^{n-1} e^{-kt} \log f_X(F^{-1}(e^{-t}))dt	      \nonumber  \\
		&=E(-\log f_X(F_X^{-1}(e^{-T_{n,k}})))	      
	\end{align}
	Ahmadi (2021) showed that the equality of the Kerridge measure of inaccuracy associated with distributions of upper and lower k-records with the underlying distribution is a characteristic of symmetric distribution. 
	In the following examples we obtain useful expressions of $H(f_{L_{n,k}},f_X)$  and $H(f_{U_{n,k}},f_X)$ in some lifetime distributions. When we take $k=1$ in Equation \ref{1}, we get the same equation obtained by Tahmasebi and Daneshi (2017).
	\begin{example}
		Let $X$ be an exponential  random variable  with rate $\theta, \ \theta>0,$ then  $H(f_{U_{n,k}},f_X)=\frac{n}{k}-\log(\theta).$ Note that for a fixed value of $n$ and $k,$ $H(f_{U_{n,k}},f_X)$ is a decreasing function of $\theta$. Also, for a fixed value of $\theta$ and $k,$ inaccuracy increases in $n.$ Also, for a fixed value of $\theta$ and $n$, inaccuracy decrease in $k$.
	\end{example}
	\begin{example}
		Let $X$ be a random variable having the Pareto distribution with pdf
		\begin{equation*}
			f_X(x)=\theta x^{-(\theta+1)}, \ x>1, \theta>0,
		\end{equation*}
		then $H(f_{U_{n,k}},f_X)=(1+\frac{1}{\theta})\frac{n}{k}-\log(\theta).$ Note that for a fixed value of $n$ and $k,$ $H(f_{U_{n,k}},f_X)$ is a decreasing function of $\theta$. Also, for a fixed value of $\theta$ and $k,$ inaccuracy increases in $n.$ Also, for a fixed value of $\theta$ and $n$, inaccuracy decrease in $k$.
	\end{example}
	
	\begin{example}
		Let X be a random variable having the Weibull distribution with pdf
		\begin{equation*}
			f_X(x)=\lambda \beta x^{\beta-1} e^{-\lambda x^{\beta}}, \ x>0, \beta>0, \lambda>0,
		\end{equation*}
		then $H(f_{U_{n,k}},f_X)=n-\log\beta-\frac{\log \lambda}{\beta}-\left(\frac{\beta-1}{\beta}\right) \frac{k^n}{\Gamma(n)} \int_{0}^{\infty} t^{n-1} e^{-kt} \log(t)dt.$
	\end{example}
	
	\begin{example}
		Let X be a random variable having a standard uniform distribution with pdf 
		\begin{equation*}
			f_X(x)=1, \ 0<x<1.
		\end{equation*}
		then $H(f_{U_{n,k}},f_X)=0 $ and  $H(f_{L_{n,k}},f_X)=0.$ 
	\end{example}
	
	\begin{example}
		Let X be a random variable having pdf 
		\begin{equation*}
			f_X(x)=3(1-x)^2, \ 0<x<1.
		\end{equation*}
		then $H(f_{U_{n,k}},f_X)=-\log3 +\frac{2n}{3k}.$ 
	\end{example}

	\begin{theorem}
		If $f(x)$ is increasing in $x,$ then $H(f_{U_{n,k}},f_X)$ is decreasing in n.
	\end{theorem}
	\textbf{Proof} Since $\frac{f_{T_{n+1,k}}(t)}{f_{T_{n,k}}(t)}=k\frac{t}{n}$ is an increasing function of $t,$ therefore  $T_{n,k}\leq_{lr}T_{n+1,k}$ and it implies $T_{n,k}\leq_{st}T_{n+1,k}$  (see, Shaked and Shanthikumar (2007)). Also, since  $f_X(F_X^{-1}(1-e^{-t}))$ is an increasing in $t$, and hence, 
	$H(f_{U_{n,k}},f_X) \geq H(f_{U_{n+1,k}},f_X).$
	
	\begin{theorem}
		If $f(x)$ is increasing in $x,$ then $H(f_{L_{n,k}},f_X)$ is increasing in $n.$
	\end{theorem}
	\textbf{Proof} Since $\frac{f_{T_{n+1,k}}(t)}{f_{T_{n,k}}(t)}=k\frac{t}{n}$ is an increasing function of $t,$ therefore  $T_{n,k}\leq_{lr}T_{n+1,k}$ and it implies $T_{n,k}\leq_{st}T_{n+1,k}$  (see, Shaked and Shanthikumar (2007)). Also, since  $f_X(F_X^{-1}(e^{-t}))$ is an decreasing in $t$, and hence, 
	$H(f_{L_{n,k}},f_X) \leq H(f_{L_{n+1,k}},f_X).$

	\section{A cumulative residual inaccuracy measure of record values}\label{section3}
	
	The residual inaccuracy measure associated with $n$th upper $k$-record value $U_{n,k}$  and parent density function $f(x)$ is given as
	\begin{align}
		H(\bar{F}_{U_{n,k}},\bar{F}_X)&=- \int_{-\infty}^{\infty} \bar{F}_{U_{n,k}}(x)\log \bar{F}_X(x)dx \nonumber  	\\
		&=  - \int_{-\infty}^{\infty}    \sum_{i=0}^{n-1} \frac{(-k \log\bar{F}_X (x) )^i}{i!} \bar{F}^k_X (x) \log \bar{F}_X(x)dx \nonumber  \\	
		&= \sum_{i=0}^{n-1} \frac{k^i}{i!}  \int_{-\infty}^{\infty}     \bar{F}^k_X (x) \left(- \log\bar{F}_X (x) \right)^{i+1} dx \nonumber \\
		&=\frac{1}{k^2} \sum_{i=0}^{n-1} (i+1) E_{U_{i+2,k}} \left(\frac{1}{\lambda_F(X)}\right) \label{HbarFUnkFx}
	\end{align}
	where $U_{i+2,k}$ is a random variable with density function \[f_{U_{i+2,k}}(x)= \dfrac{k^{i+2}}{(i+1)!}[-\log \bar{F}_X(x)]^{i+1} [\bar{F}_X(x)]^{k-1}f_X(x),\  -\infty<x<\infty \] and $\lambda_F(X)=\frac{f_X(x)}{\bar{F}_X(x)}$ is the hazard rate function of $X.$
	\begin{example}
		Let X be a random variable having a standard uniform distribution with pdf 
		\begin{equation*}
			f_X(x)=1, \ 0<x<1.
		\end{equation*}
		then 	$H(\bar{F}_{U_{n,k}},\bar{F}_X)=\sum_{i=0}^{n-1} \frac{(i+1)k^i}{(k+1)^{i+2}} .$
	\end{example}

	\begin{example}
		Let X be a random variable having an exponential distribution with pdf 
		\begin{equation*}
			f_X(x)=\theta e^{-\theta x}, \ x>0,\ \theta>0.
		\end{equation*}
		then 	$H(\bar{F}_{U_{n,k}},\bar{F}_X)=\frac{n(n+1)}{2\lambda k^2}.$
	\end{example}
	
	\begin{proposition}
		Let $X$ be an absolutely continuous non-negative random variable with survival function $\bar{F}_X$. Then, we have
		\[H(\bar{F}_{U_{n,k}},\bar{F}_X)= \sum_{i=0}^{n-1} \frac{i+1}{k}[\mu_{i+2,k}-\mu_{i+1,k}] \]
		
		where $\mu_{n,k}=\int_{0}^{+\infty} \bar{F}_{U_{n,k}}(x)dx.$
	\end{proposition}
	\textbf{Proof} Using (\ref{FUnku}), we have
	
	\[\bar{F}_{U_{i+2,k}}(x)-\bar{F}_{U_{i+1,k}}(x)=\bar{F}^k_X(x) \frac{(-k \log\bar{F}_X (x) )^{i+1}}{(i+1)!}.\]
	
	From (\ref{HbarFUnkFx}), we can write 
	\begin{align*}
		H(\bar{F}_{U_{n,k}},\bar{F}_X)&=\sum_{i=0}^{n-1} \int_{0}^{\infty} \frac{i+1}{k} \left(\bar{F}_{U_{i+2,k}}(x)-\bar{F}_{U_{i+1,k}}(x) \right) dx \\
		&=\sum_{i=0}^{n-1} \frac{i+1}{k}[\mu_{i+2,k}-\mu_{i+1,k}]
	\end{align*}
	
	\begin{proposition}
		Let $a > 0.$ For $n = 1, 2, \cdots$ it holds that 
		\[H(\bar{F}_{aU_{n,k}+b},\bar{F}_{aX+b})=a 	H(\bar{F}_{U_{n,k}},\bar{F}_X).\]
	\end{proposition}
	\textbf{Proof} Since 
	$\bar{F}_{aX+b}(x)=\bar{F}_X(\frac{x-b}{a}),$ therefore
	\begin{align*}
		H(\bar{F}_{aU_{n,k}+b},\bar{F}_{aX+b}) &=- \int_{-\infty}^{\infty} \bar{F}_{aU_{n,k}+b} (x) \log \bar{F}_{aX+b}(x)dx \\
		&= - \int_{-\infty}^{\infty} \bar{F}_{U_{n,k}} (\frac{x-b}{a}) \log \bar{F}_{X}(\frac{x-b}{a})dx \\
		&= - \int_{-\infty}^{\infty} \bar{F}_{U_{n,k}} (\frac{x-b}{a}) \log \bar{F}_{X}(\frac{x-b}{a})dx \\
		&=  -a \int_{-\infty}^{\infty} \bar{F}_{U_{n,k}} (t) \log \bar{F}_{X}(t)dt \\
		&= a H(\bar{F}_{U_{n,k}},\bar{F}_X).
	\end{align*}
	Hence, result is proved.

	\begin{proposition}
		Let $X$ be an absolutely continuous non-negative random variable with survival function $\bar{F}_X.$ Then, we have
		
		\[H(\bar{F}_{U_{n,k}},\bar{F}_X)=  \sum_{i=0}^{n-1} \frac{1}{i!} \int_{0}^{\infty} \lambda_F(z) \left[ \int_{t}^{\infty} [-k\log \bar{F}_X(x)]^i \bar{F}^k_X (x)dx \right]  dt                     \]
	\end{proposition}
	\textbf{Proof} Using equation (\ref{HbarFUnkFx}), the result $-\log\bar{F}_X(x)=\int_{0}^{x} \lambda_F(t)dt,$ and change of order of integration, we have
	\begin{align*}
		H(\bar{F}_{U_{n,k}},\bar{F}_X)&=- \int_{0}^{\infty} \bar{F}_{U_{n,k}}(x)\log \bar{F}_X(x)dx \\
		&= - \int_{0}^{\infty}    \sum_{i=0}^{n-1} \frac{(-k \log\bar{F}_X (x) )^i}{i!} \bar{F}^k_X (x) \log \bar{F}_X(x)dx \\
		&=\sum_{i=0}^{n-1} \frac{1}{i!} \int_{0}^{\infty}  \left(\int_{0}^{x} \lambda_F(t)dt \right)  (-k \log\bar{F}_X (x) )^i \bar{F}^k_X (x) dx \\
		&= \sum_{i=0}^{n-1} \frac{1}{i!} \int_{0}^{\infty} \lambda_F(z) \left[ \int_{t}^{\infty} [-k\log \bar{F}_X(x)]^i \bar{F}^k_X (x)dx \right]  dt.
	\end{align*}
	
	\begin{proposition}
		Let $X$ be an absolutely continuous non negative random variable with survival function $\bar{F}_X.$ Then, we have
		
		\[H(\bar{F}_{U_{n,k}},\bar{F}_X)=\sum_{i=0}^{n-1} \frac{1}{i!} \int_{0}^{\infty} \bar{F}^{k-1}_X (t) f_X(t) \left[\int_{0}^{t} \left(-k\log \bar{F}_X(x)\right)^{i+1}\right] dt.               \]
	\end{proposition}
	\textbf{Proof} Using (\ref{HbarFUnkFx}) and changing order of integration,  we get
	\begin{align*}
		H(\bar{F}_{U_{n,k}},\bar{F}_X)&=- \int_{0}^{\infty} \bar{F}_{U_{n,k}}(x)\log \bar{F}_X(x)dx \\
		&= - \int_{0}^{\infty}    \sum_{i=0}^{n-1} \frac{(-k \log\bar{F}_X (x) )^i}{i!} \bar{F}^k_X (x) \log \bar{F}_X(x)dx\\
		&=\frac{1}{k} \sum_{i=0}^{n-1} \frac{1}{i!} \int_{0}^{\infty}    {(-k \log\bar{F}_X (x) )^{i+1}} \bar{F}^k_X (x)dx\\
		&=\frac{1}{k} \sum_{i=0}^{n-1} \frac{1}{i!} \int_{0}^{\infty}    {(-k \log\bar{F}_X (x) )^{i+1}} \left(\int_{x}^{\infty} k \bar{F}^{k-1}_X(t)f_X(t)dt \right)dx\\
		&= \sum_{i=0}^{n-1} \frac{1}{i!} \int_{0}^{\infty}  \bar{F}^{k-1}_X(t)f_X(t)   \left( \int_{0}^{t}  {(-k \log\bar{F}_X (x) )^{i+1}} dx \right) dt.
	\end{align*}
	
	\section{A Cumulative past inaccuracy measure of record values}\label{section4}
	The Cumulative past inaccuracy measure between $X_{L_{n,k}}$
	and $X$ is presented as
	\begin{align*}
		H({F}_{L_{n,k}},{F}_X)&=- \int_{-\infty}^{\infty} F_{L_{n,k}}\log F_X(x)dx \\
		&=- \int_{-\infty}^{\infty}  \sum_{i=0}^{n-1} \frac{(-k \log F_X (x) )^i}{i!}){F}^k_X (x)  \log F_X(x)dx \\
		&= \sum_{i=0}^{n-1} \frac{k^i}{i!}  \int_{-\infty}^{\infty}     {F}^k_X (x) \left(- \log{F}_X (x) \right)^{i+1} dx \nonumber \\
		&= \sum_{i=0}^{n-1} \frac{(i+1)}{k^2} E_{L_{i+2,k}}\left[ \frac{1}{\Tilde{\lambda}(X)}\right]
	\end{align*}
	where $\Tilde{\lambda}(X)=\frac{f_X(x)}{F_X(x)}$ and $L_{i+2,k}$ is a random variable with pdf $$f_{L_{i+2,k}}(x)=\frac{k^{i+2}[-log(F_X(x)]^{i+1} [F_X(x)]^{k-1} f_X(x)}{(i+1)!}.$$
	\begin{remark}
		For $k=1,$ above expression is obtained by Tahmasebi et al. (2017).
	\end{remark}
	\begin{example}
		Let X be a random variable having a standard uniform distribution with pdf 
		\begin{equation*}
			f_X(x)=1, \ 0<x<1.
		\end{equation*}
		then $	H({F}_{L_{n,k}},{F}_X)=\frac{(i+1)k^i}{(k+1)^{i+2}} .$ 
	\end{example}
	In the following theorem, we obtain an expression of Cumulative past inaccuracy measure in terms of cdf of record values.
	\begin{theorem}
		Suppose that $X$ is a non-negative random variable with cdf $F,$ then we have
		\begin{align*}
			H({F}_{L_{n,k}},{F}_X)&= \int_{0}^{\infty} \sum_{i=0}^{n-1} \frac{(i+1)}{k} (F_{L_{i+2,k}}(x)-F_{L_{i+1,k}}(x))dx.
		\end{align*}
	\end{theorem}
	\textbf{Proof} Using (\ref{FLnku}) we obtain
	\begin{align*}
		H({F}_{L_{n,k}},{F}_X)&=- \int_{0}^{\infty} F_{L_{n,k}}\log F_X(x)dx \\
		&=- \int_{0}^{\infty}  \sum_{i=0}^{n-1} \frac{(-k \log F_X (x) )^i}{i!}){F}^k_X (x)  \log F_X(x)dx \\
		&=        \sum_{i=0}^{n-1} \frac{i+1}{k}  \int_{0}^{\infty}  \frac{(-k \log F_X (x) )^{i+1})}{(i+1)!}{F}^k_X (x)  dx\\
		&=        \sum_{i=0}^{n-1} \frac{i+1}{k}   \int_{0}^{\infty} [F_{L_{i+2,k}}(x)-F_{L_{i+1,k}}(x)]dx
	\end{align*}

	\section{Conclusion} \label{s10conclusion}
	This paper studied the inaccuracy measures between the distributions of record values and the parent random variable, particularly focusing on the Kerridge inaccuracy and cumulative residual inaccuracy measures. Based on these measures, we provided characterizations of symmetric distributions and explored their properties under various lifetime distributions. Future research can extend these ideas to other statistical records and generalized inaccuracy measures. \\

	\noindent \textbf{\Large Funding} \\
	\\
	No funding recieved.\\
	\\
	\textbf{ \Large Conflict of interest} \\
	\\
	The authors declare no conflict of interest.\\	
	\\

\end{document}